\documentclass[a4paper,10pt,conference]{ieeeconf}
\IEEEoverridecommandlockouts
\usepackage[usenames,dvipsnames]{color}
\usepackage{balance}
\usepackage{ifthen}
\usepackage{graphics}
\usepackage{subfigure}
\usepackage{graphicx,epstopdf}
\usepackage{cite}
\usepackage{amsmath,amssymb}
\usepackage{algorithm}
\usepackage{algorithmic}
\usepackage{multirow}
\usepackage{comment}

\newtheorem{theorem}{Theorem}
\newtheorem{lemma}{Lemma}

\newtheorem{corollary}{Corollary}
\newtheorem{definition}{Definition}

\newcommand{\feas}{\mathbb{F}}
\newcommand{\her}{\mathbb{H}}

\newcommand{\hN}{\mathcal{N}}
\newcommand{\hS}{\mathcal{S}}
\newcommand{\hE}{\mathcal{E}}

\newcommand{\hP}{\mathcal{P}}
\newcommand{\ph}{{\Phi}}
\newcommand{\rank}{{\mathrm{rank}}}
\newcommand{\trace}{\mathrm{tr}}
\newcommand{\re}{{\mathrm{Re}}}
\newcommand{\diag}{\mathrm{diag}}
\newcommand{\im}{{\mathrm{Im}}}
\newcommand{\range}{{\mathrm{range}}}

\newcommand{\ii}{\textbf{i}}
\newcommand{\vct}[1]{#1}
\newcommand{\mtx}[1]{#1}
\newcommand{\bq}{\begin{eqnarray}}
\newcommand{\eq}{\end{eqnarray}}
\newcommand{\bqn}{\begin{eqnarray*}}
\newcommand{\eqn}{\end{eqnarray*}}
\newcommand{\bee}{\begin{enumerate}}
\newcommand{\eee}{\end{enumerate}}
\newcommand{\bi}{\begin{itemize}}
\newcommand{\ei}{\end{itemize}}

\newboolean{showcomments}
\setboolean{showcomments}{true}
\newcommand{\lgan}[1]{\ifthenelse{\boolean{showcomments}}
{ \textcolor{red}{(Lingwen says: #1)} } {} }
\newcommand{\slow}[1]{\ifthenelse{\boolean{showcomments}}
{ \textcolor{blue}{(Steven says:  #1)}}{}}
\newcommand{\addcite}[0]{\ifthenelse{\boolean{showcomments}}
{ \textcolor{Red}{(addcite)}}{}}

\begin{document}

\title{Convex Relaxations and Linear Approximation for\\
Optimal Power Flow in Multiphase Radial Networks}

\author{Lingwen Gan and Steven H. Low
\thanks{This work was supported by NSF NetSE grant CNS 0911041, ARPA-E grant DE-AR0000226, Southern California Edison, National Science Council of Taiwan, R.O.C, grant NSC 101-3113-P-008-001, and Caltech's Resnick Institute.}
\thanks{Lingwen Gan and Steven H. Low are with the Engineering and Applied Science Division, California Institute of Technology, Pasadena, CA 91125 USA (e-mail: lgan,slow@caltech.edu).}
}	
\maketitle

\begin{abstract}
Distribution networks are usually multiphase and radial. To facilitate power flow computation and optimization, two semidefinite programming (SDP) relaxations of the optimal power flow problem and a linear approximation of the power flow are proposed. We prove that the first SDP relaxation is exact if and only if the second one is exact. Case studies show that the second SDP relaxation is numerically exact and that the linear approximation obtains voltages within 0.0016 per unit of their true values for the IEEE 13, 34, 37, 123-bus networks and a real-world 2065-bus network.
\end{abstract}


\section{Introduction}

\PARstart{A}{s} distributed generation (e.g., rooftop photovoltaic panels) and 
controllable loads (e.g.,  electric vehicles) continue to proliferate, power distribution 
systems need to become more intelligent and dynamic, better monitored and optimized.
Most distribution systems are radial, i.e., have a tree topology.
The purpose of this paper is to propose a set of models
that simplify the computation and optimization of power flows in unbalanced 
multiphase radial networks.

The optimal power flow (OPF) problem is nonconvex, and approximations and relaxations have been developed to solve it; see recent surveys in \cite{Huneault91,Momoh99a,Pandya08, Frank2012a, OPF-FERC-1, OPF-FERC-4}.
For convex relaxations, it is first proposed in \cite{Jabr2006} to solve OPF as a second-order cone programming for single-phase radial networks and in \cite{bai2008} as a semidefinite programming (SDP) for single-phase mesh networks.
 While numerically illustrated in \cite{Jabr2006} and \cite{bai2008}, whether or when the convex relaxations are exact is not
studied until \cite{Lavaei2012}; see \cite{Low2014a,Low2014b} for a survey and references
to a growing literature on  convex relaxations of OPF.

Most of these works assume a single-phase network, while distribution networks are typically multiphase and unbalanced  \cite{Kersting2006}. 
It has been observed in \cite{Berg1967, Chen1991}
that a multiphase network has an equivalent single-phase circuit
model where each bus-phase pair in the multiphase network is identified with a single
bus in the equivalent model.    
Hence methods for single-phase networks can be applied to the equivalent
model of a multiphase unbalanced network.  This approach is taken in \cite{dall2012distributed}
for solving  optimal power flow problems.
Besides, \cite{dall2012distributed} develops distributed solutions.

This paper develops convex relaxations of OPF and a linear approximation of power flow.

There are three questions on convex relaxations:
1) When can a globally optimum of OPF be obtained by solving its convex relaxation?
2) How to compute convex relaxations efficiently?
3) How to attain numerical stability?

To address 2), relaxation BIM-SDP is proposed in Section \ref{sec: BIM SDP} to improve the computational efficiency of a standard SDP relaxation by exploiting the radial network topology. While the standard SDP relaxation declares $O(n^2)$ variables where $n$ is the number of lines in the network, BIM-SDP only declares $O(n)$ variables and is therefore more efficient.

To address 3), relaxation BFM-SDP is proposed in Section \ref{sec: BFM SDP} to improve the numerical stability of BIM-SDP by avoiding ill-conditioned operations. BIM-SDP is ill-conditioned due to subtractions of voltages that are close in value. Using alternative variables, BFM-SDP avoids these subtractions and is therefore numerically more stable.

To partially address 1), we prove that BIM-SDP is exact if and only if BFM-SDP is exact, and empirically show that BFM-SDP is numerically exact for the IEEE 13, 34, 37, 123-bus networks and a real-world 2065-bus network in Section \ref{sec: case study}. Remarkably, BIM-SDP is numerically exact only for the IEEE 13 and 37-bus networks. This highlights the numerical stability of BFM-SDP.

Approximation LPF is proposed in Section \ref{sec: approximation} to estimate the voltages and power flows. LPF is accurate when line loss is small compared with power flow and voltages are nearly balanced, i.e., the voltages of different phases have similar magnitudes and differ in angle by $\sim\!\!120^\circ$. Empirically, it is presented in Section \ref{sec: case study} that LPF computes voltages within 0.0016 per unit of their true values for the IEEE 13, 34, 37, and 123-bus networks and a real-world 2065-bus network.
\section{Optimal Power Flow Problem}
\label{sec: opf}

This paper studies OPF in multiphase radial networks, and is applicable for demand response and volt/var control.

\subsection{A Standard Nonlinear Power Flow Model}
\label{subsec: BIM}

A distribution network is composed of buses and lines connecting these buses. It is usually multiphase and radial. There is a substation bus in the network with a fixed voltage. Index the substation bus by 0 and the other buses by $1,2,\ldots,n$. Let $\hN=\{0,1,\ldots,n\}$ denote the set of buses and define $\hN^+=\hN\backslash\{0\}$. Each line connects an ordered pair $(i,j)$ of buses where bus $i$ lies between bus 0 and bus $j$. Let $\hE$ denote the set of lines. Use $(i,j)\in\hE$ and $i\rightarrow j$ interchangeably. If $i\rightarrow j$ or $j\rightarrow i$, denote $i\sim j$.

Let $a,b,c$ denote the three phases of the network, let $\ph_i$ denote the phases of bus $i\in\hN$, and let $\ph_{ij}$ denote the phases of line $i\sim j$. For each bus $i\in\hN$, let $V_i^\phi$ denote its phase $\phi$ complex voltage for $\phi\in\ph_i$ and define $\vct{V}_i:=[V_i^\phi]_{\phi\in\ph_i}$; let $I_i^\phi$ denote its phase $\phi$ current injection for $\phi\in\ph_i$ and define $\vct{I}_i:=[I_i^\phi]_{\phi\in\ph_i}$; let $s_i^\phi$ denote its phase $\phi$ complex power injection for $\phi\in\ph_i$ and define $s_i:=[s_i^\phi]_{\phi\in\ph_i}$. For each line $i\sim j$, let $I_{ij}^\phi$ denote the phase $\phi$ current from bus $i$ to bus $j$ for $\phi\in\ph_{ij}$ and define $\vct{I}_{ij}:=[I_{ij}^\phi]_{\phi\in\ph_{ij}}$; let $\mtx{z}_{ij}$ denote the phase impedance matrix, assume it is full rank, and define $\mtx{y}_{ij}:= \mtx{z}_{ij}^{-1}$. 

    \begin{figure}[h!]
    \centering
    \includegraphics[scale=0.4]{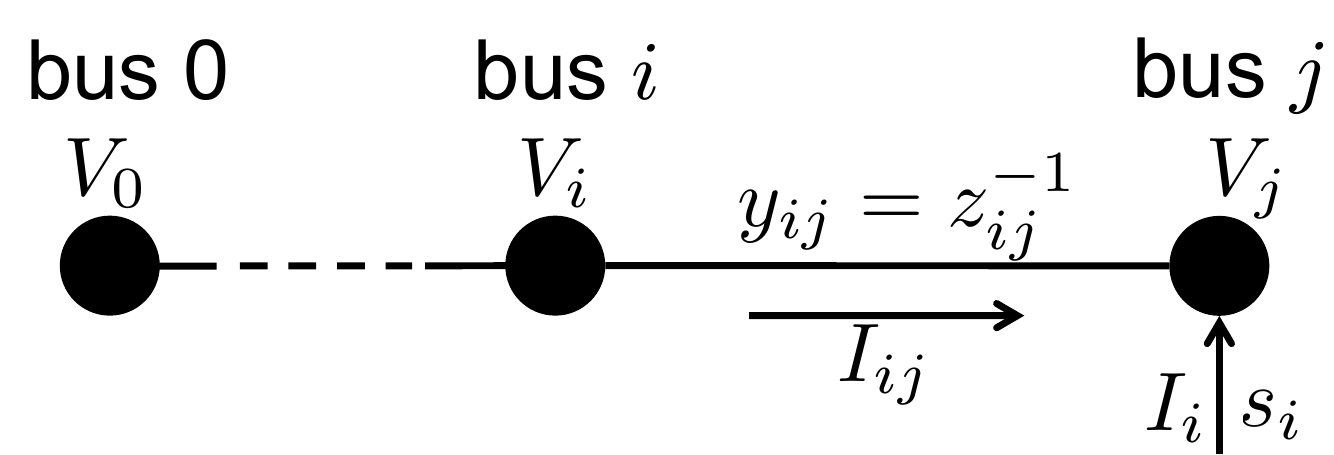}
    \caption{Summary of notations}
    \label{fig: notation}
	\end{figure}

Some notations are summarized in Fig. \ref{fig: notation}. Further, let superscripts denote projection to specified phases, e.g., if $\ph_i=abc$, then
    \[ \vct{V}_i^{ab} = (V_i^a,V_i^b)^T. \]
Fill nonexisting phase entries by 0, e.g., if $\ph_i=ab$, then
    \[ \vct{V}_i^{abc} = (V_i^a,V_i^b,0)^T. \]
Let a letter without subscripts denote a vector of the corresponding quantity, e.g., $\mtx{z}=[\mtx{z}_{ij}]_{i\sim j}$ and $s=[s_i]_{i\in\hN}$.
    
Power flows are governed by \cite{Kersting2006}:

\begin{enumerate}
\item[1)] Ohm's law:
	$\vct{I}_{ij} = \vct{y}_{ij} (\vct{V}_i^{\ph_{ij}} - \vct{V}_j^{\ph_{ij}}), \ \ i\sim j.$

\item[2)] Current balance:
	$\vct{I}_i=\sum_{j:\,i\sim j}\vct{I}_{ij}^{\ph_i}, \ \ i\in\hN.$

\item[3)] Power balance:
	$s_i=\diag(\vct{V}_i\vct{I}_i^H), \ \ i\in\hN.$
\end{enumerate}
Eliminate current variables $I_i$ and $I_{ij}$, the above model reduces to the following \emph{bus injection model} (BIM):
    \begin{equation}
    \label{BIM}
    s_i = \sum_{j:\,i\sim j} \diag\left[\vct{V}_i^{\ph_{ij}} (\vct{V}_i^{\ph_{ij}}-\vct{V}_j^{\ph_{ij}})^H\mtx{y}_{ij}^H\right]^{\ph_i}\!\!\!\!\!, \quad i\in\hN.
    \end{equation}

\subsection{Optimal Power Flow}
\label{subsec: opf}
OPF determines the power injection that minimizes generation cost subject to physical and operational constraints.

Generation cost is separable. In particular, let $C_i(s_i): \mathbb{C}^{|\ph_i|} \mapsto \mathbb{R}$ denote the generation cost at bus $i\in\hN$, and
	\begin{equation*}
	\label{objective}
	C(s) = \sum_{i\in \hN}C_i(s_i)
	\end{equation*}
is the generation cost of the network.

\begin{figure}[h!]
    \centering
    \includegraphics[scale=0.36]{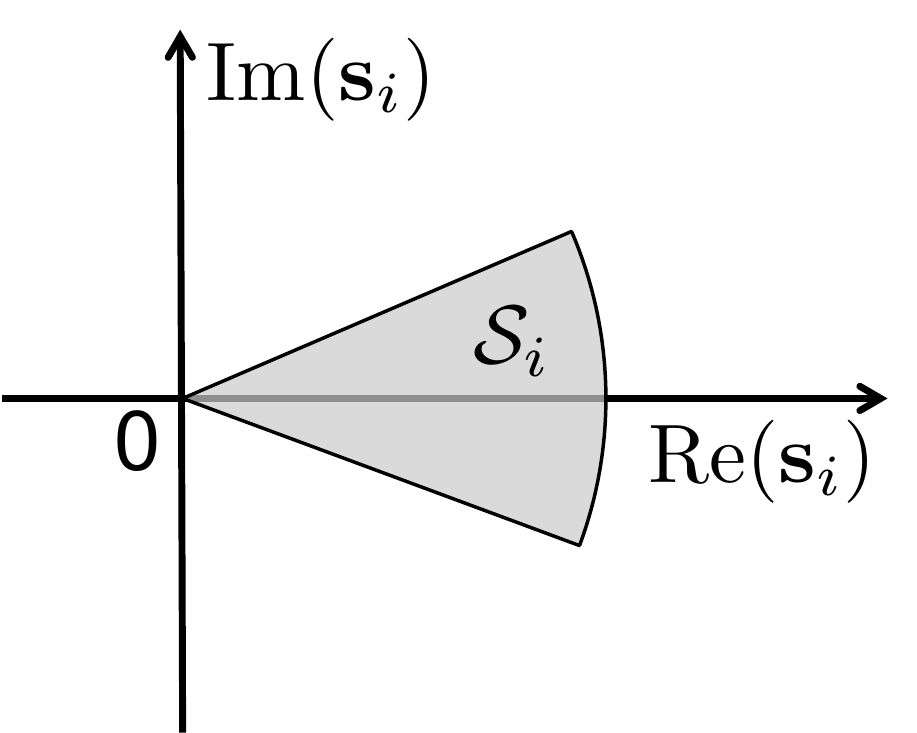}
    \hspace{0.2in}
    \includegraphics[scale=0.36]{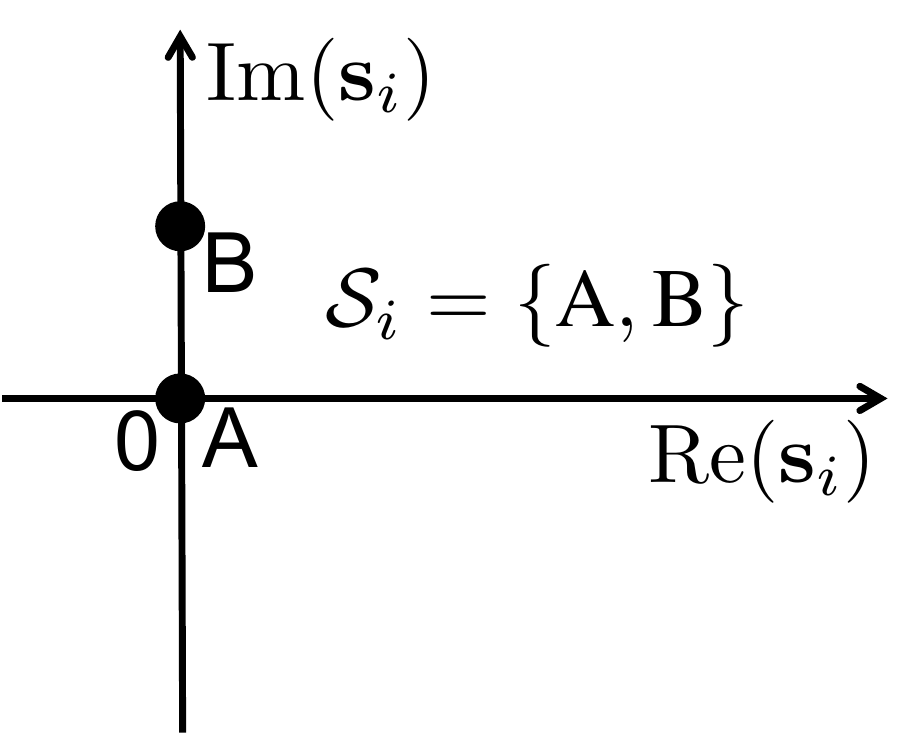}
    \caption{The left figure illustrates the set $\hS_i$ of an inverter, and the right figure illustrates the set $\hS_i$ of a shunt capacitor. Note that the set $\hS_i$ is usually not a box, and that $\hS_i$ can be nonconvex or even disconnected.}
    \label{fig: hS}
    \end{figure}

OPF has operational constraints on power injections and voltages besides physical constraints \eqref{BIM}. First, while the substation power injection $s_0$ is unconstrained, a branch bus power injection $s_i$ can only vary within some externally specified set $\hS_i$, i.e.,
	\begin{equation}\label{constraint s}
	s_i \in \hS_i, \qquad i\in \hN^+.
	\end{equation}
For example, the set $\hS_i$ of two types of devices are illustrated in Fig. \ref{fig: hS}. Note that $\hS_i$ is usually not a box, and that $\hS_i$ can be nonconvex or even disconnected.

Second, while the substation voltage $\vct{V}_0$ is fixed and given (denote by $\vct{V}_0^{\text{ref}}$ that is nonzero componentwise), a branch bus voltage can be regulated within a range, i.e., there exists $[\underline{V}_i^\phi,\overline{V}_i^\phi]_{i\in\hN^+,\phi\in\ph_i}$ such that
    \begin{subequations}
    \label{constraint v}
    \begin{align}
    & \vct{V}_0 = \vct{V}_0^{\mathrm{ref}}; \label{constraint V0}\\
	& \underline{V}_i^\phi \leq |V_i^\phi| \leq \overline{V}_i^\phi, \quad i\in \hN^+,~\phi\in\ph_i.
    \end{align}
	\end{subequations}
For example, if voltages must stay within 5\% from their nominal values, then
$0.95\leq |V_i^\phi| \leq 1.05$ per unit.

To summarize, OPF can be formulated as
    \begin{align*}
    \textbf{OPF: }\min~~ & \sum_{i\in \hN} C_i(s_i) \\
	\mathrm{over}~~ & s, \vct{V} \nonumber\\
	\mathrm{s.t.}~~ & \eqref{BIM} - \eqref{constraint v}.
    \end{align*}
The following assumptions are made throughout this paper.
\begin{enumerate}
\item The network $(\hN,\hE)$ is connected.
\item Voltage lower bounds are strictly positive, i.e.,
	\[ \underline{V}_i^\phi>0, \quad i\in\hN^+, ~\phi\in\ph_i. \]
\item Bus and line phases satisfy	
	\[ \ph_i \supseteq \ph_{ij} = \ph_j, \quad i\rightarrow j. \]
\end{enumerate}
\section{BIM Semidefinite Programming}
\label{sec: BIM SDP}
OPF is nonconvex due to \eqref{BIM}, and a standard semidefinite programming (SDP) relaxation has been developed to solve it \cite{dall2012distributed}. In this section we propose a different SDP relaxation, called BIM-SDP, that exploits the radial network topology to reduce the computational complexity of the standard SDP.

\subsection{Bus Injection Model Semidefinite Programming}
BIM-SDP is derived by shifting the nonconvexity from \eqref{BIM} to some rank constraints and removing the rank constraints.

Let $|A|$ denote the number of elements in a set $A$, and $\her^{k \times k}$ denote the set of $k\times k$ complex Hermitian matrices.
Let $v_i\in\her^{|\ph_i|\times|\ph_i|}$ for $i\in\hN$ and $W_{ij}\in\mathbb{C}^{|\ph_{ij}|\times|\ph_{ij}|}$ for $i\sim j$. If these matrices satisfy
	\[
        	\begin{bmatrix} v_i^{\ph_{ij}} & W_{ij} \\ W_{ji} & v_j \end{bmatrix}
        	= \begin{bmatrix} \vct{V}_i^{\ph_{ij}} \\ \vct{V}_j \end{bmatrix}
        	\begin{bmatrix} \vct{V}_i^{\ph_{ij}} \\ \vct{V}_j \end{bmatrix}^H,
	\quad i\rightarrow j,
    	\]
then \eqref{BIM} is equivalent to
    \[ \vct{s}_i = \sum_{j:\,i\sim j} \diag\left[ (v_i^{\ph_{ij}}-W_{ij})\mtx{y}_{ij}^H\right]^{\ph_i}, \quad i\in\hN. \]

\begin{lemma}
\label{lemma: V <=> W}
Let $v_i\in\her^{|\ph_i|\times|\ph_i|}$ for $i\in\hN$ and $W_{ij}\in\mathbb{C}^{|\ph_{ij}|\times|\ph_{ij}|}$ for $i\sim j$. If
\begin{itemize}
\item $v_0 = \vct{V}_0^\mathrm{ref} [\vct{V}_0^\mathrm{ref}]^H$ for some $\vct{V}_0^\mathrm{ref}\in\mathbb{C}^{|\ph_0|}$;
\item $\diag(v_i)$ is nonzero componentwise for $i\in\hN$;
\item $W_{ji}=W_{ij}^H$ for $i\rightarrow j$;
\item $\begin{bmatrix}
        v_i^{\ph_{ij}} & W_{ij} \\
        W_{ji} & v_j
        \end{bmatrix}$ is rank one for $i\rightarrow j$,
\end{itemize}
then Algorithm \ref{algorithm: recover V} computes the unique $\vct{V}$ that satisfies $\vct{V}_0=\vct{V}_0^\mathrm{ref}$ and
    \begin{subequations}
    \label{def v W}
    \begin{align}
    v_i &= \vct{V}_i \vct{V}_i^H, &\quad i\in\hN; \label{eq: def v}\\
    W_{ij} &= \vct{V}_i^{\ph_{ij}} (\vct{V}_j^{\ph_{ij}})^H, &\quad i\sim j.
    \end{align}
    \end{subequations}
\end{lemma}
	
	\begin{algorithm}[h!]
	\caption{Recover $V$ from $(v,W)$.}
	\label{algorithm: recover V}
	\begin{algorithmic}[1]
	
	\REQUIRE $(v,W)$ that satisfies the conditions in Lemma \ref{lemma: V <=> W}.
	\ENSURE $V$.

	\STATE $V_0 \leftarrow V_0^\mathrm{ref}$;
	\STATE $\hN_\text{visit} \leftarrow\{0\}$;
	
	\WHILE{$\hN_\text{visit} \neq \hN$}
		\STATE find $i\rightarrow j$ such that $i\in\hN_\text{visit}$ and $j\notin\hN_\text{visit}$;
		\STATE compute 
		\begin{align*}
		V_j &\leftarrow \frac{1}{\trace\left(v_i^{\ph_{ij}}\right)} W_{ji}V_i^{\ph_{ij}}; \\
		\hN_\text{visit} &\leftarrow\hN_\text{visit}\cup\{j\};
		\end{align*}
	\ENDWHILE
		
	\end{algorithmic}
	\end{algorithm}
	
	Lemma \ref{lemma: V <=> W} is proved in Appendix \ref{app: V <=> W}. It implies that OPF can be equivalently formulated as BIM-OPF. Let $A\succeq0$ denote a hermitian matrix $A$ being positive semidefinte.
	
    	\begin{subequations}
    	\label{BIM OPF}
    	\textbf{BIM-OPF: }
    	\begin{align}
    	\min~~ & \sum_{i\in \hN} C_i(\vct{s}_i) \nonumber\\
	\mathrm{over}~~ & \vct{s}_i\in\mathbb{C}^{|\ph_i|} \text{ and } v_i\in\her^{|\ph_i|\times|\ph_i|} \text{ for }i\in\hN; \nonumber\\
    	& W_{ij}\in\mathbb{C}^{|\ph_{ij}|\times|\ph_{ij}|} \text{ for }i\sim j, \nonumber\\
	\mathrm{s.t.}~~ &
	\vct{s}_i = \sum_{j:\,i\sim j} \diag\left[ (v_i^{\ph_{ij}}-W_{ij})\mtx{y}_{ij}^H\right]^{\ph_i},
	\quad i\in\hN; \label{BIM s} \\
	& \vct{s}_i \in \hS_i, \quad i\in \hN^+;  \label{BIM constraint s}\\
    	& v_0 = \vct{V}_0^\text{ref} (\vct{V}_0^\text{ref})^H; \label{BIM constraint v0}\\
	& \underline{\vct{v}}_i \leq \diag(v_i) \leq \overline{\vct{v}}_i, ~\quad i\in \hN^+;  \label{BIM constraint v} \\
    	& W_{ij} = W_{ji}^H, \quad i\rightarrow j; \label{BIM hermitian}\\
	& \begin{bmatrix}
        v_i^{\ph_{ij}} & W_{ij} \\
        W_{ji} & v_j
        \end{bmatrix} \succeq 0, \quad i\rightarrow j; \label{BIM relax}\\
	& \rank\begin{bmatrix}
        v_i^{\ph_{ij}} & W_{ij} \\
        W_{ji} & v_j
        \end{bmatrix}=1, \quad i\rightarrow j \label{BIM rank}
    \end{align}
    \end{subequations}
where the vectors $\underline{v}_i$ and $\overline{v}_i$ in \eqref{BIM constraint v} are defined as
	\[ \underline{\vct{v}}_i := [ (\underline{V}_i^\phi)^2 ]_{\phi\in\ph_i}, ~
	\overline{\vct{v}}_i := [ (\overline{V}_i^\phi)^2 ]_{\phi\in\ph_i}, \quad
	i\in\hN^+. \]
	
If $C_i$ (in the objective) and $\hS_i$ [in \eqref{BIM rank}] are convex, then BIM-OPF is convex except for \eqref{BIM rank},
and an SDP relaxation can be obtained by removing \eqref{BIM rank} from BIM-OPF.
    	\begin{align*}
    	\textbf{BIM-SDP: }\min~~ & \sum_{i\in \hN} C_i(\vct{s}_i) \\
	\mathrm{over}~~ & \vct{s}, v, W \nonumber\\
	\mathrm{s.t.}~~ & \eqref{BIM s}-\eqref{BIM relax}.  \nonumber
	\end{align*}
Note that BIM-SDP may be nonconvex due to $C_i$ and $\hS_i$.

If an optimal BIM-SDP solution $(\vct{s},v,W)$ satisfies \eqref{BIM rank}, then $(\vct{s},v,W)$ also solves BIM-OPF. Furthermore, a global optimum $(\vct{s},\vct{V})$ of OPF can be recovered via Algorithm \ref{algorithm: recover V}.

\begin{theorem}
\label{thm: BIM recover}
Given an optimal solution $(\vct{s}, v,W)$ of BIM-SDP that satisfies \eqref{BIM rank}, Algorithm 1 computes a $V$ such that $(s, V)$ solves OPF.
\end{theorem}

Theorem \ref{thm: BIM recover} follows directly from Lemma \ref{lemma: V <=> W}.
	\begin{definition}
	\label{def: BIM exact}
	BIM-SDP is \emph{exact} if every optimal solution of BIM-SDP satisfies \eqref{BIM rank}.
	\end{definition}
	
If BIM-SDP is exact, then a global optimum of OPF can be obtained by solving BIM-SDP according to Theorem \ref{thm: BIM recover}.

\subsection{Comparison with a Standard SDP}
A standard SDP relaxation of OPF has been proposed in the literature \cite{dall2012distributed}. It is derived by introducing
	\[ \tilde{W} = \begin{bmatrix}
	\vct{V}_0 \\ \vdots \\ \vct{V}_n
	\end{bmatrix}
	\begin{bmatrix}
	\vct{V}_1^H & \cdots & \vct{V}_n^H
	\end{bmatrix} \]
to shift the nonconvexity from \eqref{BIM} in BIM-OPF to
	$\rank \tilde{W} = 1,$
and removing the rank constraint. We call this relaxation standard-SDP for ease of reference.

BIM-SDP is computationally more efficient than standard-SDP since it has fewer variables. It is straightforward to verify that there are $O(n)$ variables in BIM-SDP and $O(n^2)$ variables in standard-SDP.

Standard-SDP does not exploit the radial network topology. In $\tilde{W}$,  only blocks corresponding to lines $i\sim j$ appear in other constraints than $\tilde{W}\succeq0$, i.e., if bus $i$ and bus $j$ are not connected, then block $(i,j)$ in $\tilde{W}$ only appears in $\tilde{W}\succeq0$. Since the network is radial, $n^2$ out of the $(n+1)^2$ blocks in $\tilde{W}$ only appear in $\tilde{W}\succeq0$, leaving significant potential for exploring sparsity.

Call these $n^2$ blocks that only appear in $\tilde{W}\succeq0$ the $\tilde{W}$-only blocks and the other $2n+1$ blocks the key-blocks. The role of having $\tilde{W}$-only blocks in the optimization is to make sure that the partial matrix specified by key-blocks can be completed to a positive semidefinite full matrix.

There is a standard way of positive semidefinite matrix completion---looking at the \emph{chordal} extension of a partial matrix \cite{fukuda2001exploiting}. Essentially, BIM-SDP applies this technique to exploit the radial network topology.
\section{BFM Semidefinite Programming}
\label{sec: BFM SDP}

BIM-SDP is not numerically stable and therefore a different SDP relaxation is proposed in this section.

\subsection{Alternative Power Flow Model}
We start with introducing a novel branch flow model (BFM) of power flow.
BFM enhances the numerical stability of BIM \eqref{BIM}.
BIM \eqref{BIM} is ill-conditioned due to subtractions of $V_i^{|\ph_{ij}|}$ and $V_j^{|\ph_{ij}|}$ that are close in value.
BFM attains numerical stability by avoiding such subtractions.
    
BFM is given by the following three equations.
\begin{enumerate}
\item Ohm's law:
    \begin{equation}
    \label{ohm}
    V_i^{\ph_{ij}} - V_j = z_{ij} I_{ij}, \quad i\rightarrow j.
    \end{equation}
    
\item Definition of slack variables:
\begin{equation}
    \label{def ell S}
    \ell_{ij} = I_{ij}I_{ij}^H, ~
    S_{ij} = V_i^{\ph_{ij}} I_{ij}^H, \quad
    i\rightarrow j.
    \end{equation}
    
\item Power balance:
    	\begin{align}
    	& \sum_{i:\,i\rightarrow j} \diag(S_{ij}-z_{ij}\ell_{ij}) + s_j \nonumber\\
	&\qquad\qquad\qquad = \sum_{k:\, j\rightarrow k} \diag \left( S_{jk} \right)^{\ph_j},  \quad j\in\hN. \label{power balance}
    	\end{align}
\end{enumerate}

To interpret $\ell$ and $S$, note that $\diag(\ell_{ij})$ denotes the magnitude squares of current $I_{ij}$, and $\diag(S_{ij})$ denotes the sending-end power flow on line $i\rightarrow j$. To interpret \eqref{power balance}, note that the receiving-end power flow on line $i\rightarrow j$ is
    \begin{align*}
    \diag(V_jI_{ij}^H) 
    = \diag(S_{ij} - z_{ij}\ell_{ij}).
    \end{align*}

BIM and BFM are equivalent in the sense that they share the same solution set $(\vct{s},V)$. More specifically, let
    \begin{align*}
    \feas_\text{BIM} &~:=~ \{(s,V) ~|~ (\vct{s},V) \text{ satisfies }\eqref{BIM} \}, \\
    \feas_\text{BFM} &~:=~ \left\{(\vct{s},V) ~\left|~ \begin{array} {l}
    \exists ~(I, \ell, S) \text{ such that }\\
    (s,V,I, \ell, S) \text{ satisfies \eqref{ohm}--\eqref{power balance}}
    \end{array} \right. \right\}
    \end{align*}
denote the sets of $(\vct{s},V)$ that satisfy BIM or BFM.
\begin{theorem}
\label{thm: BIM=BFM}
The solution set $\feas_\text{BIM} = \feas_\text{BFM}$.
\end{theorem}
Theorem \ref{thm: BIM=BFM} is proved in Appendix \ref{app: BIM=BFM}. It implies that OPF can be equivalently formulated as follows.
\begin{align*}
    \textbf{OPF': }\min~~ & \sum_{i\in \hN} C_i(\vct{s}_i) \\
	\mathrm{over}~~ & \vct{s}, V, I, \ell, S \nonumber\\
	\mathrm{s.t.}~~ & \eqref{constraint s} - \eqref{constraint v}, ~\eqref{ohm}-\eqref{power balance}.
    \end{align*}

\subsection{Branch Flow Model Semidefinite Programming}
A numerically stable SDP that has a similar computational efficiency as BIM-SDP is proposed in this section.

To motivate the SDP, assume \eqref{eq: def v}, \eqref{ohm}, and \eqref{def ell S} hold,
then
	\[ V_j =  V_i^{\ph_{ij}} - z_{ij} I_{ij}, \quad i\rightarrow j. \]
Multiply both sides by their Hermitian transposes to obtain
    	\begin{equation}
    	\label{v}
    	\mtx{v}_j = \mtx{v}_i^{\ph_{ij}} - ( S_{ij} z_{ij}^H + z_{ij} S_{ij}^H )
				+ z_{ij} \ell_{ij} z_{ij}^H, \quad i\rightarrow j.
	\end{equation}
Furthermore, the matrix
	\[
        \begin{bmatrix}
        \mtx{v}_i^{\ph_{ij}} & S_{ij} \\
        S_{ij}^H & \ell_{ij}
        \end{bmatrix}
        = \begin{bmatrix} V_i^{\ph_{ij}} \\ I_{ij} \end{bmatrix}
        \begin{bmatrix} V_i^{\ph_{ij}} \\ I_{ij} \end{bmatrix}^H
    \]
is positive semidefinite and rank one for $i\rightarrow j$.

\begin{lemma}
\label{lemma: v S ell <=> V}
Let $\mtx{v}_i\in\her^{|\ph_i|\times|\ph_i|}$ for $i\in\hN$. Let $S_{ij}\in\mathbb{C}^{|\ph_{ij}|\times|\ph_{ij}|}$ and $\ell_{ij}\in\her^{|\ph_{ij}|\times|\ph_{ij}|}$ for $i\rightarrow j$. If
	\begin{itemize}
	\item $\mtx{v}_0 = V_0^\mathrm{ref} [V_0^\mathrm{ref}]^H$ for some $V_0^\mathrm{ref}\in\mathbb{C}^{|\ph_0|}$;
	\item $\diag(v_i)$ is nonzero componentwise for $i\in\hN$;
	\item $(\mtx{v},S,\ell)$ satisfies \eqref{v};
	\item $\begin{bmatrix}
        		\mtx{v}_i^{\ph_{ij}} & S_{ij} \\
        		S_{ij}^H & \ell_{ij}
        		\end{bmatrix}$ is rank one for $i\rightarrow j$,
	\end{itemize}
then Algorithm \ref{algorithm: recover V I} computes the unique $(V,I)$ that satisfies $V_0=V_0^\mathrm{ref}$, \eqref{eq: def v}, \eqref{ohm}, and \eqref{def ell S}.
\end{lemma}

	\begin{algorithm}[h!]
	\caption{Recover $(V,I)$ from $(v,S,\ell)$.}
	\label{algorithm: recover V I}
	\begin{algorithmic}[1]
	
	\REQUIRE $(v,S,\ell)$ that satisfies the conditions in Lemma \ref{lemma: v S ell <=> V}.
	\ENSURE $(V,I)$.

	\STATE $V_0 \leftarrow V_0^\mathrm{ref}$;
	\STATE $\hN_\text{visit} \leftarrow\{0\}$;
	
	\WHILE{$\hN_\text{visit} \neq \hN$}
		\STATE find $i\rightarrow j$ such that $i\in\hN_\text{visit}$ and $j\notin\hN_\text{visit}$;
		\STATE compute 
		\begin{align*}
		I_{ij} &\leftarrow \frac{1}{\trace\left(v_i^{\ph_{ij}}\right)} S_{ij}^HV_i^{\ph_{ij}}; \\
		V_j &\leftarrow V_i^{\ph_{ij}} - z_{ij} I_{ij}; \\
		\hN_\text{visit} &\leftarrow\hN_\text{visit}\cup\{j\};
		\end{align*}
	\ENDWHILE
		
	\end{algorithmic}
	\end{algorithm}

Lemma \ref{lemma: v S ell <=> V} is proved in Appendix \ref{app: v S ell <=> V}. It implies that OPF' can be equivalently formulated as BFM-OPF.

    	\begin{subequations}
    	\textbf{BFM-OPF:}
    	\begin{align}
    	\min~~ & \sum_{i\in \hN} C_i(\vct{s}_i) \nonumber\\
	\mathrm{over}~~ & \vct{s}_i\in\mathbb{C}^{|\ph_i|}, \mtx{v}_i\in\her^{|\ph_i|\times|\ph_i|} \text{ for }i\in\hN; \nonumber\\
    	& S_{ij}\in\mathbb{C}^{|\ph_{ij}|\times|\ph_{ij}|},\ell_{ij}\in\her^{|\ph_{ij}|\times|\ph_{ij}|}\text{ for }i\rightarrow j, \nonumber\\
	\mathrm{s.t.}~~ & \!\!\!\sum_{i:\,i\rightarrow j}\diag(S_{ij}-z_{ij}\ell_{ij}) + s_j = \!\!\!\sum_{k:\, j\rightarrow k}\!\!\! \diag \left( S_{jk} \right)^{\ph_j}, \nonumber \\
	&\qquad\qquad\qquad\qquad   j\in\hN; \label{BFM s}\\
	& \vct{s}_i \in \hS_i, \qquad i\in \hN^+; \label{BFM constraint s} \\
    	& \mtx{v}_0 = V_0^\text{ref} (V_0^\text{ref})^H; \label{BFM constraint v0} \\
	& \underline{V}_i \leq \diag(\mtx{v}_i) \leq \overline{V}_i, \quad i\in \hN^+; \label{BFM constraint v}  \\
    	& \mtx{v}_j = \mtx{v}_i^{\ph_{ij}} - ( S_{ij} z_{ij}^H + z_{ij} S_{ij}^H ) \nonumber \\
	& \qquad\qquad + z_{ij} \ell_{ij} z_{ij}^H,\quad i\rightarrow j; \label{BFM v} \\
	& \begin{bmatrix}
        	\mtx{v}_i^{\ph_{ij}} & S_{ij} \\
        	S_{ij}^H & \ell_{ij}
        	\end{bmatrix} \succeq 0, \quad i\rightarrow j; \label{BFM relax}\\
	& \rank\begin{bmatrix}
        	\mtx{v}_i^{\ph_{ij}} & S_{ij} \\
        	S_{ij}^H & \ell_{ij}
        	\end{bmatrix}=1, \quad i\rightarrow j. \label{BFM rank}
    	\end{align}
    	\end{subequations}
	
If $C_i$ and $\hS_i$ are convex, then BFM-OPF is convex except for \eqref{BFM rank}, and an SDP relaxation can be obtained by removing \eqref{BFM rank} from BFM-OPF.
	\begin{align*}
	\textbf{BFM-SDP: }\min~~ & \sum_{i\in \hN} C_i(\vct{s}_i) \nonumber\\
	\mathrm{over}~~ & \vct{s}, \mtx{v}, S, \ell \nonumber\\
	\text{s.t.} ~~& \eqref{BFM s}-\eqref{BFM relax}.
	\end{align*}
Note that BFM-SDP may not nonconvex due to $C_i$ and $\hS_i$.

If an optimal BFM-SDP solution $(s,v,S,\ell)$ satisfies \eqref{BFM rank}, then $(s,v,S,\ell)$ also solves BFM-OPF. Moreover, Algorithm 2 produces a global optimum $(\vct{s},V, I, \ell, S)$ of OPF'.

\begin{theorem}
\label{thm: BFM recover}
Given an optimal solution $(\vct{s}, \mtx{v},S,\ell)$ of BFM-SDP that satisfies \eqref{BFM rank}, compute $(V,I)$ according to Algorithm 2. Then $(\vct{s}, V, I, \ell, S)$ solves OPF'.
\end{theorem}

Theorem \ref{thm: BFM recover} follows directly from Lemma \ref{lemma: v S ell <=> V}.

	\begin{definition}\label{def: BFM exact}
	BFM-SDP is \emph{exact} if every optimal solution of BFM-SDP satisfies \eqref{BFM rank}.
	\end{definition}
	
If BFM-SDP is exact, then a global optimum of OPF' can be obtained by solving BFM-SDP according to Theorem \ref{thm: BFM recover}.

\subsection{Comparison with BIM-SDP}
BFM-SDP is numerically more stability than BIM-SDP since it avoids subtractions of $\mtx{v}_i^{\ph_{ij}}$ and $\mtx{W}_{ij}$ that are close in value. Meanwhile, BFM-SDP has similar computational efficiency as BIM-SDP since they have the same number of variables and constraints.

There exists a bijective map between the feasible sets of BIM-SDP and BFM-SDP that preserves the objective value. Let $\feas_\text{BIM-SDP}$ and $\feas_\text{BFM-SDP}$ denote the feasible sets of BIM-SDP and BFM-SDP respectively.

\begin{theorem}
\label{thm: bijective map}
The map $f: \feas_\text{BIM-SDP} \mapsto \feas_\text{BFM-SDP}$ defined by $f(\vct{s},\mtx{v},\mtx{W})=(\vct{s},\mtx{v},S,\ell)$ where
	\begin{align*}
	S_{ij} &= (\mtx{v}_i^{\ph_{ij}} - \mtx{W}_{ij}) \mtx{y}_{ij}^H, & i\rightarrow j; \\
	\ell_{ij} &= \mtx{y}_{ij}(\mtx{v}_i^{\ph_{ij}} - \mtx{W}_{ji} - \mtx{W}_{ij} + \mtx{v}_j)\mtx{y}_{ij}^H, & i\rightarrow j
	\end{align*}
is bijective, and its inverse $g:\feas_\text{BFM-SDP} \mapsto \feas_\text{BIM-SDP}$ is given by $g(\vct{s},\mtx{v},S,\ell)=(\vct{s},\mtx{v},\mtx{W})$ where
	\begin{align*}
	\mtx{W}_{ij} = \mtx{v}_i^{\ph_{ij}} - S_{ij} z_{ij}^H, ~
	\mtx{W}_{ji} = \mtx{W}_{ij}^H, \quad
	i\rightarrow j.
	\end{align*}
\end{theorem}

Theorem \ref{thm: bijective map} is proved in Appendix \ref{app: bijective map}. It implies that $f$ is also bijective from the optimal solutions of BIM-SDP to the optimal solutions of BFM-SDP. 

\begin{corollary}
\label{cor: optimal}
Let $f$ be as in Theorem \ref{thm: bijective map}. A point $(\vct{s},\mtx{v},\mtx{W})$ solves BIM-SDP if and only if  $f(\vct{s},\mtx{v},\mtx{W})$ solves BFM-SDP.
\end{corollary}

\begin{theorem}
\label{thm: simultaneous rank one}
Let $f$ be as in Theorem \ref{thm: bijective map}. A feasible solution $(\vct{s},\mtx{v},\mtx{W})$ of BIM-SDP satisfies \eqref{BIM rank} if and only if the feasible solution $f(\vct{s},\mtx{v},\mtx{W})$ of BFM-SDP satisfies \eqref{BFM rank}.
\end{theorem}

Theorem \ref{thm: simultaneous rank one} is proved in Appendix \ref{app: simultaneous rank one}. It implies that BIM-SDP is exact if and only if BFM-SDP is exact.
\begin{corollary}
\label{cor: equivalence}
BIM-SDP is exact if and only if BFM-SDP is exact.
\end{corollary}
\section{Linear approximation}
\label{sec: approximation}
A linear approximation of the power flow LPF is proposed in this section.

LPF is obtained by assuming:
\begin{itemize}
\item[B1 ] Line losses are small, i.e.,
	$z_{ij}\ell_{ij} \ll S_{ij}$ for $i\rightarrow j$.
\item[B2 ] Voltages are nearly balanced, e.g., if $\ph_i=abc$, then
	\[ \frac{V_i^a}{V_i^b} \approx \frac{V_i^b}{V_i^c} \approx \frac{V_i^c}{V_i^a} \approx e^{j2\pi/3}. \]
\end{itemize}
With B1, omit the $z_{ij}\ell_{ij}$ terms in \eqref{BFM s} and \eqref{BFM v} to obtain
    \begin{subequations}
    \label{step1}
    \begin{align}
    & \!\!\!\sum_{i:\,i\rightarrow j}\!\!\!\diag(S_{ij}) + s_j = \!\!\! \sum_{k:\, j\rightarrow k}\!\!\! \diag ( S_{jk} )^{\ph_j}, &  j\in\hN; \label{step1 a}\\
    & v_j = v_i^{\ph_{ij}} - ( S_{ij} z_{ij}^H + z_{ij} S_{ij}^H ), & i\rightarrow j. \label{step1 b}
    \end{align}
    \end{subequations}
Given $s_j$ for $j\in\hN^+$, \eqref{step1 a} determines uniquely $s_0$ and $\diag(S_{ij})$ for $i\rightarrow j$, but not the off-diagonal entries of $S_{ij}$.
B2 is used to approximate the off-diagonal entries in $S_{ij}$ with $\diag(S_{ij})$.
Specifically, define
	\[ \alpha := e^{-j2\pi/3}, \quad
	\beta := \begin{bmatrix}
	1 \\ \alpha \\ \alpha^2
	\end{bmatrix}, \quad
	\gamma:= \begin{bmatrix}
	1 		& \alpha^2 	& \alpha \\
	\alpha 	& 1			& \alpha^2 \\
	\alpha^2	& \alpha 		& 1
	\end{bmatrix}, \]
and assume the voltages to be balanced, then
	\[ S_{ij} = V_i^{\ph_{ij}} I_{ij}^H
	\in \range(\beta^{\ph_{ij}}), \quad i\rightarrow j. \]
It follows that if $\Lambda_{ij}=\diag(S_{ij})$, let $\diag(\Lambda_{ij})$ denote a diagonal matrix with diagonal $\Lambda_{ij}$, then
	\[ S_{ij} = \gamma^{\ph_{ij}} \diag(\Lambda_{ij}). \]
	
To summarize, \eqref{step1} can be approximated by
    	\begin{subequations}
    	\label{linear approximation}
    	\begin{align}
	\textbf{LPF:}~~
    	& \sum_{i:\, i\rightarrow j} \Lambda_{ij} + s_j = \sum_{k:\, j\rightarrow k} \Lambda_{jk}^{\ph_j}, & j\in\hN; \label{step2 a}\\
	& S_{ij} = \gamma^{\ph_{ij}} \diag(\Lambda_{ij}), & i\rightarrow j; \\
    	& v_j = v_i^{\ph_{ij}} - S_{ij} z_{ij}^H - z_{ij} S_{ij}^H, & i\rightarrow j. \label{step2 c}
    	\end{align}
    	\end{subequations}
Given $s_j$ for $j\in\hN^+$ and $v_0$, \eqref{linear approximation} determines uniquely $s_0$, $(\Lambda_{ij},S_{ij})$ for $i\rightarrow j$, and $v_j$ for $j\in\hN^+$ as
	\begin{align*}
	& s_0 = - \sum_{k\in\hN^+} s_k^{\ph_0}; \\
	& \Lambda_{ij} = - \sum_{k\in\mathrm{Down}(j)} s_k^{\ph_{ij}}, & i\rightarrow j; \\
	& S_{ij} = \gamma^{\ph_{ij}} \diag(\Lambda_{ij}), & i\rightarrow j; \\
	& v_j = v_0^{\ph_j} - \!\!\! \sum_{(k,l)\in\hP_j} \left[ S_{kl} z_{kl}^H+ z_{kl} S_{kl}^H \right]^{\ph_j}, & j\in\hN^+
	\end{align*}
where $\hP_j$ denotes the path from bus 0 to bus $j$ and $\mathrm{Down}(j)$ denotes the downstream of $j$ for $j\in\hN^+$.
	
Case studies in Section \ref{sec: case study} show that LPF provides a good estimate of power flows $\Lambda$ and voltages $v$.
    
LPF generalizes the \emph{Simplified DistFlow Equations} \cite{Baran1989c} from single-phase networks to multiphase networks.   While DC approximation assumes a constant 
voltage magnitude, ignores reactive power, and assumes $r_{ij}=0$, 
LPF does not.
\section{Case studies}\label{sec: case study}
In this section, we
1) check if BIM-SDP (BFM-SDP) can be solved by the generic solver {\it sedumi} \cite{sedumi};
2) compare the running times of BIM-SDP and BFM-SDP;
3) compute how close are the BIM-SDP (BFM-SDP) solutions to rank one;
and 4) evaluate the accuracy of LPF
for the IEEE 13, 34, 37, 123-bus networks \cite{IEEE} and a real-world 2065-bus network.

The test networks are modeled by BIM and BFM with the following simplifications: 1) transformers are modeled as lines with appropriate impedances; 2) circuit switches are modeled as open or short lines depending on the status of the switch; 3) regulators are modeled as having a fixed voltage (the same as the substation); 4) distributed load on a line is modeled as two identical loads located at two end buses of the line; and 5) line shunt is modeled using the $\pi$ model---assuming a fixed impedance load at each end of the line with the impedance being half of the line shunt \cite{Kersting2006}. The real-world network locates in a residential/commercial area in Souther  California, US. All simulations are done on a laptop with Intel Core 2 Duo CPU at 2.66GHz, 4G RAM, and MAC OS 10.9.2, MATLAB R\_2013a.

\subsection{BIM-SDP vs BFM-SDP}
\label{subsec: BIM vs BFM}
OPF is set up as follows. The objective is power loss, i.e.,
	\[ C(s) = \sum_{i\in\hN} \sum_{\phi\in\ph_i} \re(s_i^{\phi}). \]
The power injection constraint \eqref{constraint s} is set up such that
\begin{enumerate}
\item for a bus $i$ representing a shunt capacitor with nameplate capacity $\overline{q}_i$,
	\[ \hS_i = \{s\in\mathbb{C}^{|\ph_i|} \mid \re(s_i)=0, ~ 0\leq\im(s_i)\leq\overline{q}_i \}; \]
\item for a solar photovoltaic bus $i$ with real power generation $p_i$ and nameplate rating $\overline{s}_i$,
	\[ \hS_i = \{s\in\mathbb{C}^{|\ph_i|} \mid \re(s_i)=p_i, ~|s_i|\leq\overline{s}_i \}; \]
\item for a bus $i$ with multiple devices, $\hS_i$ is the summation of above mentioned sets.
\end{enumerate}
Two choices of the voltage constraint \eqref{constraint v} are considered:
\begin{enumerate}
\item $\underline{V}_i^{\ph_i}=0.95$ and $\overline{V}_i^{\ph_i}=1.05$ for $i\in\hN^+$ and $\phi\in\ph_i$;
\item $\underline{V}_i^{\ph_i}=0.90$ and $\overline{V}_i^{\ph_i}=1.10$ for $i\in\hN^+$ and $\phi\in\ph_i$.
\end{enumerate}

BIM-SDP and BFM-SDP are applied to solve OPF. In particular,  the generic optimization solver \emph{sedumi} is used to solve them and results are summarized in Table \ref{table: 0.95} and \ref{table: 0.9}.

    	\begin{table}[h!]
	\caption{Simulation results using $\underline{V}=0.95$, $\overline{V}=1.05$.}
    	\label{table: 0.95}
	\scriptsize
	\vspace{-0.2in}
	\begin{center}
	\begin{tabular}{|c|c|c|c|c|c|c|}\hline
	\multirow{2}{*}{network} & \multicolumn{3}{|c|}{BIM-SDP} & \multicolumn{3}{c|}{BFM-SDP} \\
    	\cline{2-7}
    	& value & time & ratio & value & time & ratio \\ \hline
	IEEE 13-bus & 152.7 & 1.08 & 9.5e-9 & 152.7 & 0.79 & 1.6e-10 \\
	IEEE 34-bus & -100.0 & 1.97 & 1.0 & 5.001e-5 & 3.00 & 0.712 \\
	IEEE 37-bus & 212.3 & 2.32 & 1.1e-8 & 212.3 & 2.00 & 9.0e-11 \\
	IEEE 123-bus & -7140 & 6.02 & 2.2e-2 & 229.8 & 7.55 & 0.5e-11 \\
    	Rossi 2065-bus & -100.0 & 111.56 & 1.0 & 19.15 & 90.32 & 4.8e-8  \\
    	\hline
	\end{tabular}
	\end{center}
	\end{table}
	
	\begin{table}[h!]
	\caption{Simulation results using $\underline{V}=0.90$, $\overline{V}=1.10$.}
    	\label{table: 0.9}
	\scriptsize
	\vspace{-0.2in}
	\begin{center}
	\begin{tabular}{|c|c|c|c|c|c|c|}\hline
	\multirow{2}{*}{network} & \multicolumn{3}{|c|}{BIM-SDP} & \multicolumn{3}{c|}{BFM-SDP} \\
    	\cline{2-7}
    	& value & time & ratio & value & time & ratio \\ \hline
	IEEE 13-bus & 152.7 & 1.05 & 8.2e-9 & 152.7 & 0.74 & 2.8e-10 \\
	IEEE 34-bus & -100.0 & 2.22 & 1.0 & 279.0 & 1.64 & 3.3e-11 \\
	IEEE 37-bus & 212.3 & 2.66 & 1.5e-8 & 212.2 & 1.95 & 1.3e-10 \\
	IEEE 123-bus & -8917 & 7.21 & 3.2e-2 & 229.8 & 8.86 & 0.6e-11 \\
    	Rossi 2065-bus & -100.0 & 115.50 & 1.0 & 19.15 & 96.98 & 4.3e-8 \\
    	\hline
	\end{tabular}
	\end{center}
	\end{table}	

Table \ref{table: 0.95} summarizes the simulation results with $\underline{V}=0.95$ and $\overline{V}=1.05$, and Table \ref{table: 0.9} summarizes the simulations results with $\underline{V}=0.9$ and $\overline{V}=1.1$. Each table contains the (value, time, ratio) triple for each of the (network, relaxation) pairs. For example, in Table \ref{table: 0.95}, the (value, time, ratio) triple for the (BIM-SDP, IEEE 13-bus) pair is (152.7, 1.08, 9.5e-9).

The entry ``value'' stands for the objective value in the unit of kW. In the above example, with 5\% voltage flexibility, the minimum power loss of the IEEE 13-bus network computed using BIM-SDP is 152.7kW, .

The entry ``time'' stands for the running time in the unit of second. In the above example, with 5\% voltage flexibility, it takes 1.05s to solve BIM-SDP for the IEEE 13-bus network.

The entry ``ratio'' quantifies how close is an SDP solution to rank one. Due to finite numerical precision, even if BIM-SDP (BFM-SDP) is exact, its numerical solution only approximately satisfies \eqref{BIM rank} [\eqref{BFM rank}], i.e., the matrices in \eqref{BIM rank} [\eqref{BFM rank}] is only approximately rank one. To quantify how close are the matrices to rank one, one can compute their largest two eigenvalues $\lambda_1,\lambda_2$ ($|\lambda_1|\geq|\lambda_2|\geq0$) and look at their ratios $|\lambda_2/\lambda_1|$. The smaller ratios, the closer are the matrices to rank one. The maximum ratio over all matrices in \eqref{BIM rank} [\eqref{BFM rank}] is the entry ``ratio''. In the above example, with 5\% voltage flexibility, the solution of BIM-SDP for the IEEE 13-bus network satisfies $|\lambda_2/\lambda_1| \leq 9.5\times 10^{-9}$ for all matrices in \eqref{BIM rank}. Hence, BIM-SDP is numerically exact.

With 10\% voltage flexibility, BFM-SDP is numerically exact for all test networks while BIM-SDP is numerically exact only for 2 test networks. This highlights that BFM-SDP is numerically more stable than BIM-SDP, since both SDPs should be exact simultaneously if there are infinite digits of precision. When voltage flexibility reduces to 5\%, the OPF for IEEE 13-bus network becomes infeasible. Consequently, BFM-SDP is not numerically exact in this case.

To summarize, BFM-SDP is numerically exact for up to 2000-bus networks when OPF is feasible, while BIM-SDP gets into numerical difficulties for as few as 34-bus networks.

\subsection{Accuracy of LPF}
Now we evaluate the accuracy of LPF \eqref{linear approximation}. In particular, given the optimal power injections computed by BFM-SDP in Section \ref{subsec: BIM vs BFM}, we use the forward backward sweep algorithm (FBS) to obtain 
the real power flows and voltage magnitudes \cite{Kersting2007}, use LPF to estimate the power flows and voltage 
magnitudes, and compare their differences. The results are summarized in Table 
\ref{table: accuracy}.

	\begin{table}[h!]
	\caption{Accuracy of LPF.}
    	\label{table: accuracy}
	\vspace{-0.2in}
	\begin{center}
	\begin{tabular}{|c|c|c|c|c|}\hline
	\multirow{2}{*}{network} & \multicolumn{2}{|c|}{time} & \multicolumn{2}{c|}{error} \\
    	\cline{2-5}
    	& FBS & LBF & $V$ (p.u.) & $S$ (\%) \\ \hline
	IEEE 13-bus & 0.11s & 0.03s & 4.5e-4 & 3.1 \\
	IEEE 34-bus & 0.16s & 0.02s & 1.0e-3 & 4.2 \\
	IEEE 37-bus & 0.12s & 0.02s & 2.0e-4 & 1.5 \\
	IEEE 123-bus & 0.37s & 0.07s & 5.5e-4 & 3.3 \\
    	Rossi 2065-bus & 4.73s & 0.98s & 1.6e-3 & 5.3 \\
    	\hline
	\end{tabular}
	\end{center}
	\end{table}
	
It can be seen that the voltages are within 0.0016 per unit and the power flows are within 5.3\% of their true values for all test networks. This highlights the accuracy of LPF \eqref{linear approximation}.
\section{Conclusions}
Two convex relaxations, BIM-SDP and BFM-SDP, have been presented to solve OPF in multiphase radial networks. BIM-SDP explores the radial network topology to improve the computational efficiency of a standard SDP relaxation, and BFM-SDP avoids ill-conditioned operations to enhance the numerical stability of BIM-SDP. We have proved that BIM-SDP is exact if and only if BFM-SDP is exact.

A linear approximation LPF has been proposed to estimate the power flows and voltages in multiphase radial networks. LPF is accurate when line loss is small and voltages are nearly balanced. Case studies show that BFM-SDP is numerically exact and LPF obtains voltages within 0.0016 per unit of their true values for the IEEE 13, 34, 37, 123-bus networks and a real-world 2065-bus network.

\bibliographystyle{IEEEtran}
\bibliography{optimal_power_flow}

\appendix

\subsection{Proof of Lemma \ref{lemma: V <=> W}}
\label{app: V <=> W}

We prove that Algorithm \ref{algorithm: recover V} computes a $V$ that satisfies $V_0 = V_0^\mathrm{ref}$ and \eqref{def v W}. The proof of uniqueness of such $V$ is straightforward and omitted for brevity.

Let $\hN^{(0)}:=\{0\}$ and $\hN^{(k)}$ denote the set $\hN_\text{visit}$ after iteration $k=1,2,\ldots, n$ of Algorithm \ref{algorithm: recover V}. Let $\hE^{(k)}:=\hE\cap\hN^{(k)}\times\hN^{(k)}$ denote the edges of the subgraph induced by $\hN^{(k)}$ for $k=0,1,\ldots,n$.

After iteration $k\geq0$, voltage $V_i$ is recovered for $i\in\hN^{(k)}$. In particular, after iteration $n$, $V_i$ is recovered for $i\in\hN^{(n)}=\hN$. Hence, it suffices to prove
	\begin{subequations}
	\label{Lemma recover V recursive}
	\begin{align}
    	& v_i = V_i V_i^H, 				& i\in\hN^{(k)}; \\
    	& W_{ij} = V_i^{\ph_{ij}} V_j^H, 		& (i,j)\in\hE^{(k)}; \\
	& W_{ji} = V_j (V_i^{\ph_{ij}})^H, 	&\quad (i,j)\in\hE^{(k)}
    	\end{align}
	\end{subequations}
for $k=0,1,\ldots,n$.

We prove \eqref{Lemma recover V recursive} by induction. When $k=0$, \eqref{Lemma recover V recursive} holds trivially. Assume that \eqref{Lemma recover V recursive} holds for $k=K$ $(0\leq K\leq n-1)$, we prove that \eqref{Lemma recover V recursive} holds for $k=K+1$ as follows.

Let $j=\hN^{(k)} \backslash \hN^{(k-1)}$. Since the network $(\hN,\hE)$ is radial, there exists a unique $i$ such that $i\rightarrow j$. Furthermore, $i\in\hN^{(k-1)}$. It suffices to prove
	\[ v_j = V_j V_j^H, \quad
	W_{ij} = V_i^{\ph_{ij}} V_j^H, \quad
	W_{ji} = V_j (V_i^{\ph_{ij}})^H. \]
	
Since the matrix
	\[\begin{bmatrix}
	v_i^{\ph_{ij}} & W_{ij} \\
	W_{ji} & v_j
	\end{bmatrix}\]
is hermitian and rank one, there exists $\alpha,\beta\in\mathbb{C}^{|\ph_{ij}|}$ such that
	\[\begin{bmatrix}
	v_i^{\ph_{ij}} & W_{ij} \\
	W_{ji} & v_j
	\end{bmatrix}
	=\eta
	\begin{bmatrix}
	\alpha \\ \beta
	\end{bmatrix}
	\begin{bmatrix}
	\alpha^H & \beta^H
	\end{bmatrix} \]
where $\eta=\pm1$. Since $v_i^{\ph_{ij}}=V_i^{\ph_{ij}}(V_i^{\ph_{ij}})^H\succeq0$ and $v_i^{\ph_{ij}}\neq0$, one has $\eta=1$ and therefore
	\[ v_i^{\ph_{ij}} = \alpha\alpha^H, \ \
	W_{ij} = \alpha\beta^H, \ \ 
	W_{ji} = \beta \alpha^H, \ \
	v_j = \beta \beta^H. \]
Furthermore, $V_i^{\ph_{ij}}=\alpha\exp(\ii\theta)$ for some $\theta\in\mathbb{R}$ since
	\[ V_i^{\ph_{ij}} (V_i^{\ph_{ij}})^H = v_i^{\ph_{ij}} = \alpha\alpha^H. \]
It follows that
	\begin{align*}
	V_j
	&~=~ \frac{1}{\trace\left(v_i^{\ph_{ij}}\right)} W_{ji}V_i^{\ph_{ij}} \\
	&~=~ \frac{1}{\trace\left( \alpha\alpha^H \right)} \beta \alpha^H \alpha\exp(\ii\theta) \\
	&~=~ \beta \exp(\ii\theta).
	\end{align*}	
Then, it is straightforward to verify that
	\begin{align*}
	V_j V_j^H &~=~ \beta \beta^H ~=~ v_j, \\
	V_i^{\ph_{ij}}V_j^H &~=~ \alpha\beta^H ~=~ W_{ij}, \\
	V_j(V_i^{\ph_{ij}})^H &~=~ \beta\alpha^H ~=~ W_{ji}.
	\end{align*}
This completes the proof that Algorithm \ref{algorithm: recover V} computes a $V$ that satisfies $V_0 = V_0^\mathrm{ref}$ and \eqref{def v W}.

\subsection{Proof of Theorem \ref{thm: BIM=BFM}}
\label{app: BIM=BFM}

First prove $\feas_\text{BIM} \subseteq \feas_\text{BFM}$. Let $(\vct{s},V)\in\feas_\text{BIM}$, want to prove $(\vct{s},V)\in \feas_\text{BFM}$. Let $I_{ij}=y_{ij}(V_i^{\ph_{ij}}-V_j^{\ph_{ij}})$ for $i\sim j$, then $(V,I)$ satisfies \eqref{ohm}. Define $(\ell,S)$ according to \eqref{def ell S}. It suffices to prove $(s,\ell,S)$ satisfies \eqref{power balance}. This is because
	\begin{align*}
	& \sum_{k:\,j\rightarrow k} \diag(S_{jk})^{\ph_{j}}
	- \sum_{i:\,i\rightarrow j} \diag\left( S_{ij} -  z_{ij}\ell_{ij} \right) \\
	& = \sum_{k:\,j\rightarrow k} \diag\left(V_j^{\ph_{jk}} I_{jk}^H \right)^{\ph_j}
	- \sum_{i:\,i\rightarrow j} \diag\left(V_j^{\ph_{ij}} I_{ij}^H \right)^{\ph_j}  \\
	& = \sum_{i:\,i\sim j} \diag\left(V_j^{\ph_{ij}} I_{ji}^H\right)^{\ph_j} \\
	&= \sum_{i:\,i\sim j} \diag\left[V_j^{\ph_{ij}} (V_j^{\ph_{ij}}-V_i^{\ph_{ij}})^Hy_{ij}^H\right]^{\ph_j} = s_j
	\end{align*}
for $j\in\hN$. This completes the proof of $\feas_\text{BIM} \subseteq \feas_\text{BFM}$.

Next prove $\feas_\text{BFM} \subseteq \feas_\text{BIM}$.
Let $(s,V)\in\feas_\text{BFM}$, want to prove $(s,V) \in \feas_\text{BIM}$. Let $(I,\ell,S)$ be such that $(s,V,I,\ell,S)$ satisfies \eqref{ohm}--\eqref{power balance}. It suffices to prove $(s,V)$ satisfies \eqref{BIM}. This is because
	\begin{align*}
	& \sum_{i:\,i\sim j} \diag\left[V_j^{\ph_{ij}} (V_j^{\ph_{ij}}-V_i^{\ph_{ij}})^Hy_{ij}^H\right]^{\ph_j} \\
	& = \sum_{k:\,j\rightarrow k} \diag\left(V_j^{\ph_{jk}} I_{jk}^H \right)^{\ph_j}
	- \sum_{i:\,i\rightarrow j} \diag\left(V_j^{\ph_{ij}} I_{ij}^H \right)^{\ph_j} \\
	& = \sum_{k:\,j\rightarrow k} \diag(S_{jk})^{\ph_{j}}
	- \sum_{i:\,i\rightarrow j} \diag\left[ (V_i^{\ph_{ij}}-z_{ij}I_{ij})I_{ij}^H\right] \\
	& = \sum_{k:\,j\rightarrow k} \diag(S_{jk})^{\ph_{j}}
	- \sum_{i:\,i\rightarrow j} \diag\left( S_{ij} -  z_{ij}\ell_{ij} \right) = s_j
	\end{align*}
for $j\in\hN$.
This completes the proof of Theorem \ref{thm: BIM=BFM}.

\subsection{Proof of Lemma \ref{lemma: v S ell <=> V}}
\label{app: v S ell <=> V}

We prove that Algorithm \ref{algorithm: recover V I} computes a $(V,I)$ that satisfies $V_0 = V_0^\mathrm{ref}$, \eqref{eq: def v}, \eqref{ohm}, and \eqref{def ell S}. The proof of uniqueness of such $(V,I)$ is straightforward and omitted for brevity.

Let $\hN^{(0)}:=\{0\}$ and $\hN^{(k)}$ denote the set $\hN_\text{visit}$ after iteration $k=1,2,\ldots,n$ of Algorithm \ref{algorithm: recover V I}. Let $\hE^{(k)}:=\hE\cap\hN^{(k)}\times\hN^{(k)}$ denote the edges of the subgraph induced by $\hN^{(k)}$ for $k=0,1,\ldots,n$.

After iteration $k\geq0$, voltage $V_i$ is recovered for $i\in\hN^{(k)}$ and current $I_{ij}$ is recovered for $(i,j)\in\hE^{(k)}$. In particular, after iteration $n$, $V_i$ is recovered for $i\in\hN^{(n)}=\hN$ and $I_{ij}$ is recovered for $(i,j)\in\hE^{(n)}=\hE$. Hence, it suffices to prove
	\begin{subequations}
	\label{recover V I recursive}
	\begin{align}
	& v_i = V_i V_i^H, 				& i\in\hN^{(k)}; \\
	& V_i^{\ph_{ij}} - V_j = z_{ij} I_{ij}, 	& (i,j)\in\hE^{(k)}; \\
	& \ell_{ij} = I_{ij}I_{ij}^H, 			& (i,j)\in\hE^{(k)}; \\
	& S_{ij} = V_i^{\ph_{ij}}I_{ij}^H, 		& (i,j) \in\hE^{(k)}
	\end{align}
	\end{subequations}
for $k=0,1,\ldots,n$.

We prove \eqref{recover V I recursive} by induction. When $k=0$, \eqref{recover V I recursive} holds trivially. Assume that \eqref{recover V I recursive} holds for $k=K$ ($0\leq K\leq n-1$), we prove that \eqref{recover V I recursive} holds for $k=K+1$ as follows.

Let $j= \hN^{(k)}\backslash\hN^{(k-1)}$. Since the network $(\hN,\hE)$ is radial, there exists a unique $i$ such that $i\rightarrow j$. Furthermore, $i\in\hN^{(k-1)}$. It suffices to prove
	\begin{align*}
	& v_j = V_j V_j^H,
	& V_i^{\ph_{ij}} - V_j = z_{ij}I_{ij}, \\
	&\ell_{ij} = I_{ij}I_{ij}^H,
	& S_{ij} = V_i^{\ph_{ij}} I_{ij}^H.
	\end{align*}
	
Since the matrix
	\[\begin{bmatrix}
	v_i^{\ph_{ij}} & S_{ij} \\
	S_{ij}^H & \ell_{ij}
	\end{bmatrix}\]
is hermitian and rank one, there exists $\alpha,\beta\in\mathbb{C}^{|\ph_{ij}|}$ such that
	\[\begin{bmatrix}
	v_i^{\ph_{ij}} & S_{ij} \\
	S_{ij}^H & \ell_{ij}
	\end{bmatrix}
	=\eta
	\begin{bmatrix}
	\alpha \\ \beta
	\end{bmatrix}
	\begin{bmatrix}
	\alpha^H & \beta^H
	\end{bmatrix}	\]
where $\eta=\pm1$. Since $v_i^{\ph_{ij}}=V_i^{\ph_{ij}}(V_i^{\ph_{ij}})^H\succeq0$ and $v_i^{\ph_{ij}}\neq0$, one has $\eta=1$ and therefore
	\[ v_i^{\ph_{ij}} = \alpha\alpha^H, \ \
	S_{ij} = \alpha\beta^H, \ \ 
	\ell_{ij} = \beta \beta^H. \]
Furthermore, $V_i^{\ph_{ij}}=\alpha\exp(\ii\theta)$ for some $\theta\in\mathbb{R}$ since
	\[ V_i^{\ph_{ij}} (V_i^{\ph_{ij}})^H = v_i^{\ph_{ij}} = \alpha\alpha^H. \]
It follows that
	\begin{align*}
	I_{ij}
	&~=~ \frac{1}{\trace\left(v_i^{\ph_{ij}}\right)} S_{ij}^HV_i^{\ph_{ij}} \\
	&~=~ \frac{1}{\trace\left( \alpha\alpha^H \right)} \beta \alpha^H \alpha\exp(\ii\theta) \\
	&~=~ \beta \exp(\ii\theta), \\
	V_j
	&~=~ V_i^{\ph_{ij}} - z_{ij}I_{ij} \\
	&~=~ \alpha\exp(\ii\theta) - z_{ij}  \beta \exp(\ii\theta) \\
	&~=~ \left( \alpha - z_{ij}  \beta \right) \exp(\ii\theta).
	\end{align*}
Then, it is straightforward to verify that
	\begin{align*}
	V_j V_j^H &~=~ ( \alpha - z_{ij}  \beta ) ( \alpha - z_{ij}  \beta )^H \\
	&~=~ v_i^{\ph_{ij}} - (z_{ij} S_{ij}^H + S_{ij}z_{ij}^H ) + z_{ij} \ell_{ij} z_{ij}^H ~=~ v_j, \\
	I_{ij}I_{ij}^H &~=~ \beta\beta^H ~=~ \ell_{ij}, \\
	V_i^{\ph_{ij}}I_{ij}^H &~=~ \alpha\beta^H ~=~ S_{ij}.
	\end{align*}
This completes the proof that Algorithm \ref{algorithm: recover V I} computes a $(V,I)$ that satisfies $V_0 = V_0^\mathrm{ref}$, \eqref{eq: def v}, \eqref{ohm}, and \eqref{def ell S}.

\subsection{Proof of Theorem \ref{thm: bijective map}}
\label{app: bijective map}

First prove that $f(s,v,W) \in \feas_\text{BFM-SDP}$ for any $(s,v,W)\in\feas_\text{BIM-SDP}$. Let $(s,v,W) \in \feas_\text{BIM-SDP}$, let $(s,v,S,\ell)=f(s,v,W)$, want to prove $(s,v,S,\ell)\in\feas_\text{BFM-SDP}$.

It is straightforward that $(s,v,S,\ell)$ satisfies \eqref{BFM constraint s}--\eqref{BFM constraint v}. The point $(\vct{s}, v, \mtx{S}, \mtx{\ell})$ satisfies \eqref{BFM s} because
	\begin{align*}
	& S_{ij} - z_{ij} \ell_{ij} \\
	&= (v_i^{\ph_{ij}} - W_{ij})y_{ij}^H - (v_i^{\ph_{ij}} - W_{ij} - W_{ji} + v_j)y_{ij}^H \\
	&= - (v_j - W_{ji}) y_{ij}^H
	\end{align*}
for $i\rightarrow j$ and therefore
	\begin{align*}
	& \sum_{k:\,j\rightarrow k} \diag(S_{jk})^{\ph_j} - \sum_{i:\,i\rightarrow j}\diag\left(S_{ij}-z_{ij}\ell_{ij} \right) \\
	&= \sum_{k:\,j\rightarrow k} \diag\left[(v_j^{\ph_{jk}}- W_{jk})y_{jk}^H\right]^{\ph_j} \\
	&\qquad + \sum_{i:\,i\rightarrow j}\diag\left[(v_j - W_{ji})y_{ij}^H\right] \\
	&= \sum_{i:\,i\sim j}\diag\left[ ( v_j^{\ph_{ij}}  -  W_{ji} )\mtx{y}_{ji}^H \right]^{\ph_j} ~=~ \vct{s}_j
	\end{align*}
for $j\in\hN$. The point $(\vct{s}, v, \mtx{S}, \mtx{\ell})$ satisfies \eqref{BFM v} because
	\begin{align*}
	& v_i^{\ph_{ij}} - (\mtx{S}_{ij}\mtx{z}_{ij}^H + \mtx{z}_{ij} \mtx{S}_{ij}^H) + \mtx{z}_{ij}\mtx{\ell}_{ij}\mtx{z}_{ij}^H \\
	&= v_i^{\ph_{ij}} - (v_i^{\ph_{ij}} - W_{ij} + v_i^{\ph_{ij}} - W_{ji}) \\
	& \qquad + v_i^{\ph_{ij}} - W_{ij} - W_{ji} + v_j \\
	&= v_j
	\end{align*}
for $i\rightarrow j$. The point $(\vct{s}, v, \mtx{S}, \mtx{\ell})$ satisfies \eqref{BFM relax} because
	\begin{align}
	& \begin{bmatrix}
	v_i^{\ph_{ij}} & \mtx{S}_{ij} \\
	\mtx{S}_{ij}^H & \mtx{\ell}_{ij}
	\end{bmatrix}
	\succeq 0 \nonumber \\
	\Leftrightarrow ~ & v_i^{\ph_{ij}} \succeq 0, ~S_{ij}\in\range(v_i^{\ph_{ij}}),~\mtx{\ell}_{ij} \succeq \mtx{S}_{ij}^H (v_i^{\ph_{ij}})^+\mtx{S}_{ij} \nonumber\\
	\Leftrightarrow ~ & v_i^{\ph_{ij}} \succeq 0, ~W_{ij}\in\range(v_i^{\ph_{ij}}),
	~v_i^{\ph_{ij}} - W_{ij} - W_{ji} + v_j \nonumber\\
	&\qquad \succeq (v_i^{\ph_{ij}} - W_{ji}) (v_i^{\ph_{ij}})^+ (v_i^{\ph_{ij}} - W_{ij}) \nonumber\\
	\Leftrightarrow ~ & v_i^{\ph_{ij}} \succeq 0, ~W_{ij}\in\range(v_i^{\ph_{ij}}), ~v_j \succeq W_{ji} (v_i^{\ph_{ij}})^+ W_{ij} \nonumber\\
	\Leftrightarrow ~ & \begin{bmatrix}
	v_i^{\ph_{ij}} & W_{ij} \\
	W_{ji} & v_j
	\end{bmatrix}
	\succeq 0 \label{PSD equivalent}
	\end{align}
for $i\rightarrow j$. This completes the proof that $f(s,v,W) \in \feas_\text{BFM-SDP}$ for any $(s,v,W)\in\feas_\text{BIM-SDP}$.

Next show that $g(s,v,S,\ell)\in\feas_\text{BIM-SDP}$ for any $(s,v,S,\ell) \in \feas_\text{BFM-SDP}$. Let $(s,v,S,\ell) \in \feas_\text{BFM-SDP}$, let $g(s,v,S,\ell) = (s,v,W)$, want to prove $(s,v,W)\in\feas_\text{BIM-SDP}$.

It is straightforward that $(s,v,W)$ satisfies \eqref{BIM constraint s}--\eqref{BIM hermitian}. The point $(\vct{s}, v, W)$ satisfies \eqref{BIM s} because
	\begin{align*}
	 (v_j - W_{ji}) y_{ij}^H &~=~ -(v_i^{\ph_{ij}} - z_{ij} S_{ij}^H - v_j ) \mtx{y}_{ij}^H \\
	&~=~ -(S_{ij}z_{ij}^H - z_{ij}\ell_{ij}z_{ij}^H ) \mtx{y}_{ij}^H \\
	&~=~ -(S_{ij}  - z_{ij}\ell_{ij})
	\end{align*}
for $i\rightarrow j$ and therefore
	\begin{align*}
	& \sum_{i:\, i\sim j} \diag\left[ (v_j^{\ph_{ij}} - W_{ji})y_{ji}^H \right]^{\ph_j} \\
	&= \sum_{i:\, i\rightarrow j} \diag\left[ (v_j - W_{ji})y_{ji}^H \right] \\
	&\qquad + \sum_{k:\, j\rightarrow k} \diag\left[ (v_j^{\ph_{jk}} - W_{jk})y_{jk}^H \right]^{\ph_j} \\
	&= - \sum_{i:\, i\rightarrow j}\!\! \diag( S_{ij} - z_{ij}\ell_{ij} ) 
	+ \sum_{k:\, j\rightarrow k}\!\! \diag(S_{jk} )^{\ph_j} = \vct{s}_j
	\end{align*}
for $j\in\hN$. 
The point $(\vct{s}, v, W)$ satisfies \eqref{BIM relax} due to \eqref{PSD equivalent}. This completes the proof that $g(s,v,S,\ell)\in\feas_\text{BIM-SDP}$ for any $(s,v,S,\ell) \in \feas_\text{BFM-SDP}$.

It is straightforward to verify that $f\circ g$ and $g\circ f$ are both identity maps. This completes the proof of Theorem \ref{thm: bijective map}.

\subsection{Proof of Theorem \ref{thm: simultaneous rank one}}
\label{app: simultaneous rank one}
Let $(\vct{s}, v, W)\in\feas_\text{BIM-SDP}$ and $(\vct{s}, v, \mtx{S}, \mtx{\ell}) = f(\vct{s}, v, W)$. It suffices to prove that
	\[
		\rank
		\begin{bmatrix}
        		v_i^{\ph_{ij}} & S_{ij} \\
        		S_{ij}^H & \ell_{ij}
        		\end{bmatrix}
		=1
		~\Longleftrightarrow~
		\rank
		\begin{bmatrix}
        		v_i^{\ph_{ij}} & W_{ij} \\
        		W_{ji} & v_j
        		\end{bmatrix}
		=1
	\]
for $i\rightarrow j$.

Fix an arbitrary $i\rightarrow j$. Since $\rank(v_i^{\ph_{ij}})\geq1$ by \eqref{BIM constraint v},
	\begin{align*}
	& \rank\begin{bmatrix}
        	v_i^{\ph_{ij}} & S_{ij} \\
        	S_{ij}^H & \ell_{ij}
        	\end{bmatrix}=1 \\
	\Leftrightarrow ~ & \rank(v_i^{\ph_{ij}})=1,
	~S_{ij}\in\range(v_i^{\ph_{ij}}),
	~\ell_{ij}=\mtx{S}_{ij}^H (v_i^{\ph_{ij}})^+\mtx{S}_{ij} \\
	\Leftrightarrow ~ & \rank(v_i^{\ph_{ij}})=1,
	~W_{ij}\in\range(v_i^{\ph_{ij}}), \\
	& v_i^{\ph_{ij}} - W_{ji} - W_{ij} + v_j \\
	& \qquad= (v_i^{\ph_{ij}} - W_{ji}) (v_i^{\ph_{ij}})^+(v_i^{\ph_{ij}} - W_{ij})\\
	\Leftrightarrow ~ & \rank(v_i^{\ph_{ij}})=1,
	~W_{ij}\in\range(v_i^{\ph_{ij}}), \\
	& v_j = W_{ji} (v_i^{\ph_{ij}})^+ W_{ij} \\
	\Leftrightarrow ~ & \rank\begin{bmatrix}
        	v_i^{\ph_{ij}} & W_{ij} \\
        	W_{ji} & v_j
        	\end{bmatrix}=1.
	\end{align*}
This completes the proof of Theorem \ref{thm: simultaneous rank one}.

\end{document}